\documentclass[11pt]{amsart} 
\usepackage{verbatim, latexsym, amssymb, amsmath,color}
\usepackage{epsfig}

\newtheorem{theorem}{Theorem}[section]
\newtheorem{corollary}[theorem]{Corollary}

\theoremstyle{definition}

\newtheorem{remark}[theorem]{Remark}

\title{Minimal surfaces - variational theory and applications}
\author{Fernando Cod\'{a} Marques}
\address{Instituto de Matem\'atica Pura e Aplicada (IMPA) \\ Estrada Dona Castorina 110 \\ 22460-320 Rio de Janeiro \\ Brazil}
\email{coda@impa.br}

\thanks{The author is grateful to  \'{E}cole Polytechnique, \'{E}cole Normale Sup\'{e}rieure, and Universit\'{e} Paris-Est (Marne-la-Vall\'{e}e)
for the hospitality during the writing of this paper.}

\begin{document}

\begin{abstract}
Minimal surfaces  are among the most natural objects in Differential Geometry, and have been studied for the past 250 years ever since the pioneering work of Lagrange. The subject is characterized by a profound beauty, but perhaps even more remarkably, minimal surfaces (or minimal submanifolds) have encountered  striking applications in other fields, like three-dimensional topology, mathematical physics, conformal geometry, among others. Even though it has been the subject of intense activity, many basic open problems still remain. In this lecture we will survey recent advances in this area and discuss some future directions. We will give special  emphasis to the variational aspects of the theory as well as to the applications to other fields.  
\end{abstract}



\maketitle

\section{Introduction and Results}

Minimal submanifolds are  solutions of the most basic variational problem of submanifold geometry, that of extremizing the area. This was first considered by Lagrange (1762), who raised the question of existence of surfaces of least area having a given closed curve in three-space as the boundary. He derived the differential
equation that must be satisfied by a function of two variables whose graph  minimizes area among surfaces with a given contour. Later Meusnier discovered that this is equivalent to the vanishing of the mean curvature, and the study of the differential geometry of these surfaces was started. The theory of minimal surfaces (or minimal submanifolds) has been developed over the years by several outstanding mathematicians, and it  is now extremely rich. It has been extended to other ambient geometries, and it is full of beautiful examples and deep  theorems. 

This paper attempts to give an overview and discuss some recent advances in the subject, with emphasis on the variational aspects and applications. Being such a large field, we do not have  the pretension of being exhaustive. Part of the material of this article is also discussed in the contribution of Andr\'{e}  Neves \cite{neves}.

Let us begin with a discussion of the first variation formula.

\medskip

\subsection{First Variation Formula}
Let $\Sigma$ be a two-dimensional oriented surface in $\mathbb{R}^3$, and let $N$ denote a unit normal field.  The local geometry of $\Sigma$ at a point $p$ can be understood in terms of the principal curvatures $k_1,k_2$, the maximum and minimum curvatures  of the intersections of the surface with  normal planes passing through $p$. The classical notions of curvature of a surface in three-space are: 
\begin{itemize}
\item the mean curvature $H=(k_1+k_2)/2$, 
\item the Gauss curvature $K=k_1\cdot k_2$. 
\end{itemize}

The Gauss curvature $K$, according to the Theorema Egregium of Gauss, is an intrinsic notion, i.e., depends only on measurements made in the surface $\Sigma$ without reference to the space in which the surface is embedded. The mean curvature $H$, by contrast, is extrinsic and is naturally related to the area functional as follows.

Given a smooth variation $F:(-\varepsilon, \varepsilon) \times \Sigma \rightarrow \mathbb{R}^3$  of $\Sigma$, with $F(0, \cdot) = {\rm id}$, $\Sigma_t=F(t,\Sigma)$, the {\it First Variation Formula} tells us that
$$
\frac{d}{dt}_{|t=0} {\rm area}(\Sigma_t)= - \int_\Sigma \langle \vec{H}, X \rangle \, d\Sigma + \int_{\partial \Sigma} \langle \nu, X\rangle \, ds,
$$
where $\vec{H}=H\cdot N$ is the mean curvature vector of $\Sigma$ in $\mathbb{R}^3$, $\nu$ is the outward unit conormal vector of $\partial \Sigma$, and $X=\frac{\partial F}{\partial t}(0,\cdot)$ is the variational vector field. 

The
formula holds true in the more general setting of a $k$-dimensional submanifold $\Sigma$ immersed in an $n$-dimensional  Riemannian manifold $M$. It leads us to define  a {\it minimal submanifold} as one  for which the mean curvature vector vanishes ($\vec{H}=0$) or, equivalently, one for which the first derivative of area is zero with respect to any variation that keeps the boundary fixed ($X=0$ on $\partial \Sigma$).

\medskip

\subsection{Plateau's Problem}
Minimal surfaces can be physically represented as
soap films, which can be experimentally produced by dipping wire contours into soapy water. These experiments were systematically carried out
by the physicist Joseph Plateau in the 19th century. The problem, raised by Lagrange, of finding the surface of least area with a given boundary in Euclidean space became known as the Plateau's Problem. 

This became a central question in the field, inspiring the development of a great amount of mathematics since the time of Riemann, until it was independently solved in 1930 by Douglas \cite{douglas} and Rad\'{o} \cite{rado}. They considered two-dimensional surfaces that were given by mappings of the unit disk, and proved in particular that every smooth Jordan curve in Euclidean space is the boundary of a least area surface of disk type. Later Morrey \cite{morrey}
extended this existence theory to two-dimensional surfaces in $n$-dimensional Riemannian manifolds that are homogeneously regular, a condition satisfied by all closed manifolds.

The search for solving the Plateau's problem in greater generality, extending it  to submanifolds  of higher dimensions and of arbitrary topological type lead to the development of Geometric Measure Theory.  In a seminal paper, Federer and Fleming \cite{federer-fleming} introduced the  class of integral currents to model $k$-dimensional domains of integration.  This class had the right compactness properties to allow the solution of an extremely general Plateau's problem by the direct method of the calculus of variations.  The regularity of these $k$-area minimizing currents was the subject of much work later on (\cite{almgren66, almgren.big, bombieri-degiorgi-giusti, degiorgi2, delellis-spadaro, federer-fleming, federer.singular, simons}). In the case of codimension one, the area minimizing current is smooth outside a singular set of codimension 7.

An important source of area minimizing submanifolds comes from the calibration theory introduced by Harvey and Lawson \cite{harvey-lawson}. Complex submanifolds in K\"{a}hler manifolds and special Lagrangian submanifolds in Calabi-Yau manifolds are calibrated, therefore area minimizing in their homology class.

Area minimizing submanifolds are in particular {\it stable}, i.e., the second variation of area is nonnegative for any
variational normal vector field $X$ with $X=0$ on $\partial \Sigma$. 


\medskip

\subsection{Minimizing in homotopy and isotopy classes}
The existence of incompressible minimal surfaces in Riemannian manifolds  was proven by Schoen and Yau \cite{Schoen-Yau2}  and independently by Sacks and Uhlenbeck \cite{sacks-uhlenbeck2}. Let $\Sigma_g$ be a compact Riemann surface of genus $g$, and suppose $f:\Sigma_g \rightarrow M$ is a continuous map such that the action $f_{*}:\pi_1(\Sigma_g)\rightarrow \pi_1(M)$  induced at the level of the fundamental groups is injective.  Then  there exists a branched minimal immersion $h:\Sigma_g \rightarrow M$ that minimizes area among all maps $h':\Sigma_g \rightarrow M$ that satisfy $h'_*=f_*$. The idea is to first fix the conformal structure of $\Sigma_g$ and minimize the Dirichlet energy
$E(f)=\int_{\Sigma_g} |df|^2 \, d\mu$. This produces a family of harmonic maps, and the second step is to minimize their energies over the Teichm\"{u}ller space.
In case ${\rm dim}(M)=3$, the works of Osserman \cite{osserman} and Gulliver \cite{gulliver} establish that the map $h$ has no branch points, i.e., it is a smooth immersion. 

In \cite{msy}, Meeks, Simon and Yau proved the existence of embedded minimal surfaces by minimizing area (instead of energy) in nontrivial isotopy classes. Their methods, together with Schoen's curvature estimates  \cite{schoen}, are used to establish regularity in some treatments of min-max theory (see \cite{colding-delellis}). The embeddedness question and applications to three-dimensional topology were also the subject of \cite{meeks-yau80, meeks-yau.plateau}.

An existence theory that produces minimal two-spheres in every compact Riemannian manifold was developed in the  article of Sacks and Uhlenbeck \cite{sacks-uhlenbeck}. These minimal two-spheres are parametrized conformally by a harmonic map from $S^2$ to $M$ that is an immersion (not necessarily an embedding) outside finitely many branch points. Every harmonic map from $S^2$ to $M$ is also a conformal branched minimal immersion, hence these minimal spheres can be constructed by extremizing the Dirichlet energy:
$
E(f)=\int_{S^2} |df|^2 \, d\mu,
$
$f:S^2\rightarrow M$. A major difficulty arises from the facts that the energy $E$ is conformally invariant and that the group of conformal transformations of $S^2$ is noncompact. These issues are dealt with in \cite{sacks-uhlenbeck}  by approximating the energy $E$ by a family of energy functionals $E_\alpha$, $\alpha>1$, that satisfy the Palais-Smale condition. The possible loss of compactness by concentration of energy in the limit, as $\alpha \rightarrow 1$, is treated by the introduction of the renormalization (or blow-up) technique. 


\medskip

\subsection{Scalar curvature and Positive Mass Conjecture} A celebrated application of minimal hypersurfaces of minimizing type to mathematical physics is the proof of the Positive Mass Conjecture by Schoen and Yau \cite{Schoen-Yau1, Schoen-Yau3}. Later Witten \cite{Witten} gave a different proof using harmonic spinors. The theorem establishes that the total mass of an isolated gravitational system, modeled by an asymptotically flat spacetime obeying the dominant energy condition, must be positive unless  the spacetime is the Minkowski space (of zero mass).  In the time-symmetric case, this reduces to showing that the mass of an asymptotically flat Riemannian three-manifold of nonnegative scalar curvature is positive unless the manifold is the Euclidean space. 

The proof of Schoen and Yau is by contradiction. If the mass is negative, they prove that one can construct (by taking a limit of solutions to Plateau problems) a complete orientable area-minimizing (hence stable) minimal surface $\Sigma$ in $M$. Curvature estimates for stable minimal submanifolds are needed in this process (\cite{schoen}, \cite{schoen-simon-yau}). By  the Second Variation Formula, the stability condition gives that
$$
\int_{\Sigma}\left(|\nabla f|^2 - (|A|^2 + Ric(N,N))f^2\right) \, d\Sigma \geq 0
$$
for any smooth function $f$ with compact support in $\Sigma$. Here $A$ denotes the second fundamental form. 
The idea is to exploit the stability inequality and arrive at a contradiction with the Gauss-Bonnet Theorem. 

The same type of argument works in the compact setting to prove that the torus $T^3$ does not admit a metric of positive scalar curvature (\cite{Schoen-Yau2}).  The trick is to use the Gauss equation to make the ambient scalar curvature appear, rewriting the stability condition as:
 $$
\int_{\Sigma}\left(|\nabla f|^2 - (\frac12 R_M - K_\Sigma+ \frac12 |A|^2)f^2\right) \, d\Sigma \geq 0
$$ 
for $f\in C^\infty_0(\Sigma)$. We denote by $R_M$ the scalar curvature of $M$. Since any Riemannian three-torus contains a stable minimal $T^2$, by minimization in a   homotopy class, we obtain a contradiction between $R_M>0$ and the Gauss-Bonnet Theorem by choosing $f\equiv 1$.

Gromov and Lawson \cite{gromov-lawson} used spinorial techniques to prove that the torus $T^n$ does not carry a metric of positive scalar curvature, for any $n$.
The argument of Schoen and Yau extends to any dimension $3\leq n\leq 7$, and breaks down in higher dimensions because of the possibility of singularities in the solution to the Plateau problem. The proof of Witten of the positive mass theorem works in any dimension $n\geq 3$ under the topological requirement that the manifold is spin.
 Despite recent approaches, the positive mass conjecture is still open for nonspin manifolds in high dimensions.

Finally, we point out that minimal surfaces also play a very important role in general relativity by modeling apparent horizons of black holes. The Penrose inequality (proven by Huisken and Ilmanen \cite{huisken-ilmanen} and Bray \cite{bray}), for instance, gives a beautiful and sharp inequality between the total mass and the area of an outermost minimal sphere.

\medskip

\subsection{Min-max methods}

So far we have discussed only minimization questions, but in general minimal submanifolds are critical points of saddle type. 
Poincar\'{e} \cite{poincare} realized the importance of constructing these critical points, in the context of geodesics. 
He  asked the foundational question of whether every Riemannian two-sphere contains a closed geodesic. Geodesics are, of course, examples of minimal submanifolds.  This question had a tremendous impact in mathematics. Firstly, it can be interpreted from two different points of view: as the search for  periodic orbits of the geodesic flow or for critical points of the length functional.  Secondly, entirely new techniques and topological ideas had to be developed to answer it.

The first breakthrough was due to Birkhoff \cite{birkhoff}, who introduced min-max methods to the problem. He defined the notion of {\it sweepout}: a continuous family of closed curves $\{c_t\}_{t \in [0,1]}$ in $S^2$ that can be written as 
$$
c_t = f\left(\{x \in S^2: x_3=1-2t\}\right)
$$
for some degree one map $f:S^2 \rightarrow S^2$. Given a Riemannian two-sphere $(S^2,g)$, one can consider the min-max invariant 
$$
L = \inf_f \sup_{t\in [0,1]} L(c_t),
$$
where $L(c)$ denotes the length of the curve $c$. Birkhoff proved that $L>0$ and that $L=L(\gamma)$ for some smooth closed geodesic $\gamma$. The geodesic $\gamma$ is obtained as a limit of curves of maximal length in a minimizing sequence of sweepouts. Therefore every Riemannian two-sphere $(S^2,g)$ contains at least one closed geodesic. The fact that every compact Riemannian manifold $M^n$ contains a closed geodesic was later established by Lyusternik and Fet \cite{lyusternik-fet}. 

The work of Birkhoff inspired the development of Morse theory  and Lusternik-Schnirelman theory, fundamental ideas in mathematics that  brought together the fields of topology and the calculus of variations. 
For instance, Lusternik and Schnirelmann \cite{lusternik}  introduced  new topological methods that lead to a proof that every metric on a two-sphere admits at least three simple (embedded) closed geodesics. The search for a rigorous proof motivated a great amount of work  (see \cite{taimanov1}). The proof of Grayson \cite{grayson} uses a parabolic partial differential equation, the curve shortening flow of curves. In the early 1990s,  this activity culminated with the proof that every Riemannian two-sphere contains infinitely many smooth closed geodesics. This follows by combining the works of Franks \cite{franks} and Bangert \cite{bangert}, while Hingston \cite{hingston} proved  quantitative results. 

It is natural to ask whether every closed Riemannian $n$-manifold $M$ contains a $k$-dimensional closed minimal submanifold. This suggests looking for a Morse theory for minimal varieties, similar in spirit to the case of closed geodesics. The first step was done in \cite{almgren}, where Almgren (by suggestion of Federer) computes the homotopy groups of the space $\mathcal{Z}_k(M)$ of $k$-dimensional integral cycles  (integral currents with boundary zero) of $M$.  Almgren proved that the $l$-dimensional homotopy group of $\mathcal{Z}_k(M)$ is isomorphic to the $(k+l)$-dimensional homology group $H_{k+l}(M,\mathbb{Z})$. This theorem gives examples of one-parameter homotopically nontrivial families of cycles, so one could think of applying min-max methods just as in the case of Birkhoff's sweepouts. One major problem is that the mass functional (area functional) is only lower semicontinuous in the flat topology (the natural topology for currents). This is not a drawback in questions of minimization, but it could lead to loss of area in the limit, thereby preventing an unstable minimal surface from being detected by the min-max approach. (A study of the existence of unstable solutions to the two-dimensional Plateau problem had been done by Morse and Tompkins in \cite{morse-tompkins}.)

In \cite{almgren-varifolds}, Almgren deals with this issue by considering the measure theoretic class of surfaces called varifolds. The natural topology of the space of varifolds allowed for both good compactness properties and continuity of the area functional. Almgren devised a very general min-max theory that worked in any dimension and codimension, and for families of cycles of any number of parameters. He was able to show that every closed Riemannian manifold $M^n$ contains at least one stationary integral $k$-dimensional varifold, for each $1\leq k\leq n$. (Recall that a varifold is stationary when the first variation of area is zero with respect to any smooth deformation of the ambient space.) This Morse-theoretic theorem left open the question of regularity of the min-max minimal variety. Note that a general stationary integral varifold contains an open dense set where it is a smooth submanifold, according to the regularity theorem of Allard \cite{allard}.

In \cite{pitts}, Pitts improved considerably the theory of Almgren by showing that the stationary integral varifold can be chosen to satisfy an additional variational property, the almost minimizing condition. Roughy speaking, an almost minimizing varifold is one that can be arbitrarily approximated by integral currents that nearly minimize area. If the codimension is one, curvature estimates for stable minimal hypersurfaces can  be used to prove regularity. Pitts employed the pointwise curvature estimates of Schoen, Simon and Yau \cite{schoen-simon-yau} to prove that if $3\leq n \leq 6$ then the min-max minimal variety can be chosen to be an embedded smooth closed minimal hypersurface. Schoen and Simon \cite{schoen-simon} then extended the regularity theory through different methods to higher dimensions, allowing singular sets of codimension $7$ (see Wickramasekera \cite{wickramasekera} for a general regularity theory of stable hypersurfaces). By combining these results, the theorem is:

\begin{theorem}\label{pitts.theorem}
Let $(M^n,g)$ be a compact Riemannian manifold, with $n\geq 3$. Then $M$ contains a stationary integral varifold $\Sigma$, whose support is smooth outside a singular set of codimension 7.  In particular, if $n\leq 7$ then $M$ contains a smooth embedded closed minimal hypersurface $\Sigma^{n-1}$.
\end{theorem}

\begin{remark}
If $M^n$ satisfies $H_{n-1}(M,\mathbb{Z})\neq 0$, then the existence of $\Sigma$ follows by direct minimization of area inside a nontrivial homology class $\sigma \in H_{n-1}(M,\mathbb{Z})$, together with the regularity theory for codimension one area minimizing currents mentioned before. 
\end{remark}

\begin{remark}
Although it is true, as mentioned before, that every compact surface contains a  closed geodesic, the Almgren-Pitts min-max theory applied to that setting does not yield a smooth object. The three-legged starfish example (see \cite{almgren-varifolds}) suggests a situation in which the min-max curve has the shape of a figure eight. 
The Almgren-Pitts theory has been applied to the setting of one-cycles by Calabi and Cao \cite{calabi-cao}, who proved that in two-spheres with nonnegative curvature the closed geodesics of shortest length are always embedded.
\end{remark}

A beautiful application of higher index minimal two-spheres was given by Micallef and Moore \cite{micallef-moore}, following  the Sacks-Uhlenbeck existence approach. They introduced the positive isotropic curvature condition for Riemannian manifolds $M^n$ of dimension $n\geq 4$, and proved that it implies $\pi_2(M)=\cdots=\pi_{[n/2]}(M)=0$ when $M$ is compact. If $M$ is further assumed to be simply connected, it follows that it must be homeomorphic to a sphere.  The method is reminiscent of Synge's theorem in Riemannian geometry. The idea is
to analyze the index of the Sacks-Uhlenbeck harmonic spheres (by applying Morse theory to the perturbed energy functionals) and compare with a lower bound for the index coming from the Riemann-Roch theorem. The trick is to use the complexified version of the second variation formula, previously considered by Siu and Yau \cite{siu-yau} in their proof of the Frankel conjecture. For a general compact Riemannian manifold of positive isotropic curvature, the conjecture is that its fundamental group must contain a free subgroup of finite index (Fraser \cite{fraser}, Gromov \cite{gromov.macroscopic}). In dimension four the conjecture follows from the complete classification of Chen, Tang and Zhu \cite{chen-tang-zhu}. We also mention that minimal surfaces have been recently used by Liu \cite{liu} to classify complete three-manifolds of nonnegative Ricci curvature, thereby completing the work of Schoen and Yau \cite{schoen-yau.ricci}.

\vspace{0.5cm}


\subsection{Recent applications}\label{recent.applications} We briefly describe some recent applications of min-max minimal surfaces that will be discussed in more detail in Sections \ref{finite.time}, \ref{willmore.mobius} and \ref{further}.

Colding and Minicozzi \cite{Colding-Minicozzi1, Colding-Minicozzi2}  found a splendid application of min-max methods to the study of the topology of compact three-manifolds. Their contribution fits together with the celebrated proof of the Poincar\'{e} conjecture by Perelman \cite{PERELMAN02,PERELMAN03A,PERELMAN03B}, achieved by a profound analysis of Hamilton's Ricci flow. In \cite{PERELMAN03B}, Perelman proved  that the Ricci flow with surgeries starting at a homotopy three-sphere becomes extinct in finite time. The finite time extinction result was actually proven more generally for any closed orientable three-manifold whose prime decomposition contains only non-aspherical factors.  In this case one can avoid the analysis of the longtime behavior of the  flow, and conclude that $M$ must be diffeomorphic to a connected sum of copies of $S^2 \times S^1$ and of spherical space forms $S^3/\Gamma$.


The proof of Perelman of the finite time extinction result was itself based in minimal surface techniques of a variational nature. For technical reasons it required in addition a regularized version of the curve shortening flow, due to  
Altschuler and Grayson \cite{altschuler-grayson}. In \cite{Colding-Minicozzi1, Colding-Minicozzi2}, Colding and Minicozzi provided an alternative elegant argument inspired by the Sacks-Uhlenbeck approach.  The evolution equation of the area of the  min-max minimal surface produced by sweeping  out a homotopy 3-sphere by $S^2$'s implies that the area eventually decreases to zero in finite time. We will say more about their argument in Section \ref{finite.time}.



Recently, the author and A. Neves have been able to find a connection between the min-max theory of minimal surfaces in $S^3$ and the  Willmore conjecture (1965). The main new insights come from the analysis of the geometric and topological properties of a canonical 5-dimensional family of surfaces in the three-sphere, to which we apply min-max theory for the area functional. 

The Willmore energy of a closed surface $\Sigma$ immersed in Euclidean three-space is the total integral of the square of the mean curvature:
$$
\mathcal{W}(\Sigma) = \int_\Sigma H^2 d\Sigma.
$$
This functional is invariant under the action of any  conformal transformation of $\mathbb{R}^3$, and it appears naturally in physical contexts. It was proposed in the 1800s by Sophie Germain \cite{germain} to describe elastic shells, and more modernly it appears in the Helfrich model \cite{helfrich} of mathematical biology as one of the terms that contribute to the energy of cell membranes.

Willmore proved that the round spheres (of energy $4\pi$) minimize the Willmore energy among all closed surfaces in $\mathbb{R}^3$, and asked what is the optimal shape among surfaces of some   fixed topological type. Motivated by the analysis of circular tori of revolution, Willmore made a conjecture for the case of genus one: 
\medskip
 
\noindent {\bf Willmore Conjecture (1965, \cite{willmore})}: \textit{The integral of the square of the mean curvature of a torus immersed in $\mathbb{R}^3$ is at least $2\pi^2$.}

\medskip

The torus  $\Sigma_{\sqrt{2}}$, obtained by rotation of a circle of radius 1 with center at distance $\sqrt{2}$ of the axis of revolution, satisfies $\mathcal{W}(\Sigma_{\sqrt{2}})=2\pi^2.$ 

The author and A. Neves proved the following theorem that implies the conjecture:
\begin{theorem}[\cite{marques-neves}]\label{theorem.A.euclidean}
Let $\Sigma\subset \mathbb{R}^3$ be a closed, embedded smooth surface with  genus $g\geq 1$. Then $\mathcal{W}(\Sigma)\geq 2\pi^2$, and $\mathcal{W}(\Sigma)= 2\pi^2$ if and only if $\Sigma$ is a conformal image of $\Sigma_{\sqrt{2}}$. \end{theorem}

We were motivated by the problem of producing the Clifford torus $\hat{\Sigma}=S^1(\frac{1}{\sqrt{2}}) \times S^1(\frac{1}{\sqrt{2}})$, a minimal surface of area $2\pi^2$ and  index five in $S^3$, by min-max methods. As a byproduct of the construction, we proved:
\begin{theorem}[\cite{marques-neves}] \label{theorem.B}
Let $\Sigma\subset S^3$ be a closed, embedded smooth surface with  genus $g\geq 1$. If $\Sigma$ is minimal then ${\rm area}(\Sigma)\geq 2\pi^2$, and ${\rm area}(\Sigma)=2\pi^2$ if and only if $\Sigma$ is the image of the Clifford torus under an ambient isometry. 
\end{theorem}

Together with I. Agol and A. Neves, we used these ideas to solve another conformally invariant variational problem.
Recall that a 2-component link in $\mathbb{R}^3$ is a pair $(\gamma_1,\gamma_2)$ of rectifiable curves $\gamma_i: S^1 \rightarrow \mathbb{R}^3$, $i=1,2$,  such that $\gamma_1(S^1) \cap \gamma_2(S^1)=\emptyset$. The {\it M\"{o}bius cross energy} of the link $(\gamma_1,\gamma_2)$ is defined to be
$$
E(\gamma_1,\gamma_2) = \int_{S^1 \times S^1} \frac{|\gamma_1'(s)||\gamma_2'(t)|}{|\gamma_1(s)-\gamma_2(t)|^2}\, ds\, dt.
$$ 
This energy was introduced in \cite{freedman-he-wang} and has the property, like the Willmore energy of surfaces, of conformal invariance.

It was conjectured by Freedman-He-Wang (1994, \cite{freedman-he-wang}) that the M\"{o}bius energy of any nontrivial link  should be at least $2\pi^2$. The equality is attained by the stereographic projection of the so-called {\it standard Hopf link}:
$$\hat{\gamma}_1(s)=(\cos s, \sin s,0,0) \in S^3\quad\mbox{and}\quad\hat{\gamma}_2(t)=(0,0,\cos t,\sin t) \in S^3.$$ 
It follows from a result of He \cite{He02} that it suffices to prove the conjecture for links $(\gamma_1,\gamma_2)$ that have linking number ${\rm lk}(\gamma_1,\gamma_2)=\pm 1$. This is what we proved in \cite{agol-marques-neves}:

\begin{theorem}[\cite{agol-marques-neves}]\label{fhw.conjecture.theorem} Let $\gamma_i: S^1 \rightarrow \mathbb{R}^3$, $i=1,2$, be a 2-component link in $\mathbb{R}^3$ with $|{\rm lk}(\gamma_1,\gamma_2)| = 1$.  Then $E(\gamma_1,\gamma_2) \geq 2\pi^2$. 

Moreover, if $E(\gamma_1,\gamma_2)=2\pi^2$ then there exists a conformal map $F:\mathbb{R}^4 \rightarrow \mathbb{R}^4$ such that $(F\circ \gamma_1,F\circ \gamma_2)$ describes the standard Hopf link up to orientation.
\end{theorem}

We will give an idea about the proofs of these statements and what they have in common in Section \ref{willmore.mobius}.

Motivated by the results in the case of geodesics, Yau conjectured in \cite{yau1} (first problem in the Minimal Surfaces section) that every compact  Riemannian three-manifold admits an infinite number of smooth, closed, immersed minimal surfaces. In \cite{marques-neves-infinitely}, we have been able to prove this conjecture in the positive Ricci curvature setting, or more generally, for manifolds $(M,g)$ that satisfy the {\it embedded Frankel property}: 
\begin{itemize}
\item any two smooth, closed, embedded minimal hypersurfaces of $M$ intersect each other. 
\end{itemize}

The author and A. Neves proved:
\begin{theorem}[\cite{marques-neves-infinitely}]\label{infinitely.theorem} Let $(M,g)$ be a compact Riemannian manifold of dimension $(n+1)$, with $2\leq n\leq 6.$ Suppose that $M$ satisfies the embedded Frankel property. Then $M$ contains an infinite number of distinct smooth, closed, embedded, minimal hypersurfaces. 
\end{theorem}

Since manifolds of positive Ricci curvature satisfy the embedded Frankel property \cite{frankel}, we derive the following corollary:

\begin{corollary}[\cite{marques-neves-infinitely}]\label{ricci.thm} Let $(M,g)$ be a compact  Riemannian $(n+1)$-manifold with $2\leq n\leq 6.$ If the Ricci curvature of $g$ is positive, then $M$ contains an infinite number of distinct smooth, closed, embedded, minimal hypersurfaces. 
\end{corollary}

The proof of Theorem \ref{infinitely.theorem} uses the Almgren-Pitts min-max theory for the area functional, combined with ideas from Lusternik-Schnirelmann theory. The idea is to  apply min-max theory to the multiparameter families of hypersurfaces (mod $2$ cycles) studied by Gromov  \cite{gromov0, gromov, gromov2}  and Guth \cite{guth}. 

In the case of generic metrics on three-manifolds, there is  an alternative approach described by  Kapouleas  \cite{kapouleas2} to construct an infinite number of  embedded minimal surfaces. This works by either desingularizing two intersecting minimal surfaces or by doubling an existing unstable minimal surface. We also point out that Kahn and Markovic \cite{kahn-markovic} proved the existence of infinitely many incompressible surfaces in compact hyperbolic three-manifolds. Hence, by minimization, every Riemannian metric on one of these manifolds admits infinitely many smooth immersed minimal surfaces.

An informal overview of the proof of Theorem \ref{infinitely.theorem} will be given in Section \ref{further}.

\medskip

\subsection{Other advances}
Many other techniques have been employed in the study of minimal surfaces, like mononicity formulas, the strong maximum principle, curvature estimates, Weierstrass representation and others. We describe  recent and important advances in minimal surface theory that make use of some of these methods. We refer the reader to the book of Meeks and P\'{e}rez \cite{meeks-perez} and the references therein for more results of that type.

Colding and Minicozzi 
have developed a  deep theory about the structure of arbitrary embedded  minimal surfaces of bounded genus. As a consequence, they showed (\cite{colding-minicozzi.calabi})  that any complete embedded minimal surface with finite topology in $\mathbb{R}^3$ must be proper, thereby proving that the Calabi-Yau conjectures for embedded surfaces are true. In particular, these results imply that the examples constructed by  Jorge and Xavier \cite{jorge-xavier} and by Nadirashvili \cite{nadirashvili} cannot be embedded.

An important theme in minimal surface theory is the search for classification results. Meeks and Rosenberg \cite{meeks-rosenberg} have used the theory of Colding and Minicozzi to prove a great classification theorem:  the plane and the helicoid are the only complete properly embedded simply connected minimal surfaces in $\mathbb{R}^3$. (The properness assumption can be removed because of \cite{colding-minicozzi.calabi}.) This generalizes the classical  Bernstein theorem: the only entire minimal graph in $\mathbb{R}^3$ is the plane. (The Bernstein theorem holds true in $\mathbb{R}^n$ for $n\leq 8$, and it is false if $n\geq 9$.) The Bernstein  theorem has another  beautiful extension, proven earlier and independently by do Carmo and Peng \cite{do-carmo-peng}, Fischer-Colbrie and Schoen \cite{fischer-colbrie-schoen}, and Pogorelov \cite{pogorelov}: the only complete orientable stable minimal surface in $\mathbb{R}^3$ is the plane. The nonorientable case has been recently settled by  Ros \cite{ros.onesided}, but  the classification of stable minimal hypersurfaces in $\mathbb{R}^n$ for $n\geq 4$ is still an open problem.

Another important recent contribution to minimal surface theory fits perfectly into the classification theme. In \cite{brendle}, Brendle proved the Lawson conjecture (\cite{lawson}): the only minimal embedded torus in $S^3$ is the Clifford torus. The proof is based on the maximum
principle technique, and it extends the work of Andrews \cite{andrews} (on mean curvature flow) in an insightful way. The technique was used again by Andrews and Li \cite{andrews-li} to confirm the classification of constant mean curvature tori in $S^3$ conjectured by Pinkall and Sterling \cite{pinkall-sterling}. 

It is also fundamental to enrich the class of known examples of minimal surfaces. In the 1980s, Costa \cite{costa} discovered a very important one: the first finite topology properly embedded minimal surface in $\mathbb{R}^3$ to be discovered after the plane, the catenoid and the helicoid. The embeddedness was proven by Hoffman and Meeks \cite{hoffman-meeks}, who constructed an infinite family of similar examples. Recently, Hoffman, Traizet and White \cite{hoffman-traizet-white} have presented a construction, by variational means, of new minimal surfaces in $\mathbb{R}^3$:  embedded helicoidal minimal surfaces of every genus $g$. These  surfaces are limits of minimal surfaces in the homogeneous spaces $S^2(r) \times \mathbb{R}$ as $r \rightarrow \infty$. There has been a lot of activity recently in the study of minimal and constant mean curvature surfaces in homogeneous three-manifolds (see Fern\'{a}ndez and Mira \cite{fernandez-mira} for a survey). Minimal surfaces in $\mathbb{H} \times \mathbb{R}$, for instance, have been used by Collin and Rosenberg \cite{collin-rosenberg} to construct harmonic diffeomorphisms from the complex plane $\mathbb{C}$ onto the hyperbolic plane $\mathbb{H}$. 

Minimal surfaces in $S^3$ with arbitrary genus were constructed by Lawson \cite{lawson70}. Other examples were found by Karcher, Pinkall and Sterling \cite{ karcher-pinkall-sterling}. This list has been recently enlarged with the examples of Kapouleas and Yang \cite{kapouleas-yang} (by gluing techniques) and of Choe and Soret \cite{choe-soret} (by Lawson's method). 
Minimal surfaces can also be constructed using degree arguments. In \cite{white.varying}, White develops a mapping degree theory for closed minimal surfaces that implies every metric of positive Ricci curvature in $S^3$ contains an embedded minimal torus and at least two embedded minimal spheres (the existence of one minimal sphere follows by min-max \cite{smyth}). 

Finally, we mention that there are several interesting connections between the theory of minimal hypersurfaces and the study of partial differential equations that come from phase transition theory. This is a very active field. We refer the reader to the article of Pacard \cite{pacard} and the references therein.

\section{Finite time extinction of Ricci flow}\label{finite.time}

 Let $M^3$ be a closed orientable three-manifold that is prime and non-aspherical. In particular, $\pi_3(M) \neq 0$ by standard topology. We consider the space $\Omega$ of continuous maps $$\varphi: S^2 \times [0,1] \rightarrow M$$ 
such that $\varphi(S^2 \times \{0\})$ and $\varphi(S^2 \times \{1\})$ are points and so that  $$t\in [0,1] \rightarrow \varphi(\cdot, t)$$ is a continuous map from $[0,1]$ to $C^0\cap W^{1,2}$. Here $W^{1,2}$ denotes the Sobolev space of maps $f:S^2 \rightarrow M$ with first derivatives in $L^2$. 

We denote by $\Omega_\varphi$ the subspace of maps $\psi\in \Omega$ that are homotopic to a given $\varphi \in \Omega$ through maps in $\Omega$. A map $\varphi \in \Omega$ naturally induces a continuous map $\hat{\varphi}:S^3 \rightarrow M$. We take $\Omega_\varphi$ so that the corresponding map $\hat{\varphi}$ represents a nontrivial element of $\pi_3(M)$. In that case we say that elements of $\Omega_\varphi$ are sweepouts of $M$.

Given a Riemannian metric $g$ on $M$, the {\it energy width} of $\Omega_\varphi$ with respect to $g$ is defined to be the min-max invariant:
$$
W_E=W_E(\varphi,g) = \inf_{\psi\in \Omega_\varphi} \sup_{t\in [0,1]} E(\psi(\cdot,t)).
$$
Recall that
$
E(f) = \frac12 \int_{S^2} |df|^2 \, dv_{\overline{g}},
$
where $\overline{g}$ is the standard metric of curvature one on $S^2$.
If $\varphi$ is a sweepout of $M$, then $W_E>0$. 
We can define the {\it area width} of $\Omega_\varphi$ with respect to $g$ similarly by
$$
W_A=W_A(\varphi,g) = \inf_{\psi\in \Omega_\varphi} \sup_{t\in [0,1]} {\rm area}(\psi(\cdot,t)),
$$
where
$
{\rm area}(f) = \int_{S^2} {\rm Jac}(f)\, dv_{\overline{g}}.
$

If $\{e_1,e_2\}$ is an orthonormal basis of $T_pS^2$, then 
$$
{\rm Jac}(f)(p)= \sqrt{|df(e_1)|^2|df(e_2)|^2-\langle df(e_1),df(e_2) \rangle^2} \leq \frac12 (|df(e_1)|^2+|df(e_2)|^2),
$$
hence ${\rm area}(f)\leq E(f)$ with equality if and only if $f$ is almost conformal, i.e., $\langle df(e_1),df(e_2) \rangle=0$ and $|df(e_1)|=|df(e_2)|$ almost everywhere. It follows immediately that $W_A \leq W_E$.  By conformally reparametrizing, through the Riemann mapping theorem for variable metrics (see \cite{jost}), the maps $\varphi_i(\cdot, t)$ of a minimizing sequence $\varphi_i\in \Omega_\varphi$ for $W_A$: 
$$
\lim_{i \rightarrow \infty} \sup_{t\in [0,1]} {\rm area}(\varphi_i(\cdot,t))=W_A,
$$
Colding and Minicozzi showed that $W_E\leq W_A$. Hence $W_E=W_A$. In particular, any sequence of sweepouts $\varphi_i\in \Omega_\varphi$ that is minimizing for the energy, i.e., such that $$\lim_{i \rightarrow \infty} \sup_{t\in [0,1]} E(\varphi_i(\cdot,t))=W_E,$$
is also minimizing for the area: $
\lim_{i \rightarrow \infty} \sup_{t\in [0,1]} {\rm area}(\varphi_i(\cdot,t))=W_A.
$

The  theorem below establishes the existence of good minimizing sequences of sweepouts. Recall that any harmonic map $u:S^2 \rightarrow M$ is a conformal branched minimal immersion.
\begin{theorem}[\cite{Colding-Minicozzi2}]\label{good.sweepouts}
Suppose that $\varphi:S^2 \times [0,1]\rightarrow M$ induces a homotopically nontrivial map $\hat{\varphi}: S^3 \rightarrow M^3$. For any Riemannian metric $g$ on $M$, there exists a sequence of sweepouts $\varphi_i\in \Omega_\varphi$ with $$\lim_{i \rightarrow \infty} \sup_{t\in [0,1]} E(\varphi_i(\cdot,t))=W_E=W$$
and such that any sufficiently large slice (in area) is close in varifold sense to a union of branched minimal spheres. More precisely, for any $\varepsilon>0$ we can find $i_0\in \mathbb{N}$ and $\delta>0$  so that if $i\geq i_0$ and $s \in [0,1]$ satisfy
$$
{\rm area}(\varphi_i(\cdot,s))\geq W-\delta,
$$
then 
$$
{\bf F}\left(\varphi_i(\cdot,s), \cup_k \{u_k\}\right)\leq \varepsilon
$$
for some  finite union $\cup\{u_k\}$ of harmonic maps
$u_k:S^2 \rightarrow M$. 
\end{theorem}
In the above statement ${\bf F}$ denotes the ${\bf F}$-metric, that gives the weak topology to the space of varifolds. We identify a map $f:S^2 \rightarrow M$ with the two-dimensional rectifiable varifold $f_\#(S^2)$ it induces in $M$.

The sequence $\varphi_i$ is obtained from an arbitrary minimizing sequence $\psi_i$ for the energy through harmonic replacement, a process similar in spirit to Birkhoff's curve shortening deformation \cite{birkhoff}.  The convergence to the branched minimal spheres is inspired by compactness theorems for harmonic maps with bounded energy (\cite{jost},  \cite{parker}, \cite{parker-wolfson}, \cite{siu-yau}).  The fact that the convergence is in varifold sense implies, in particular, that there is no loss of area or energy in the limit.

Suppose now that we evolve the metric on $M^3$ by Hamilton's Ricci flow \cite{hamilton}:
$$
\frac{\partial g}{\partial t}=-2 Ric_{g(t)}.
$$
Let $[0,T)$ be the maximal time of smooth existence. If $\varphi:S^2 \times [0,1]\rightarrow M$ induces a homotopically nontrivial map $\hat{\varphi}: S^3 \rightarrow M^3$, the function $W(t)=W_E(\varphi,g(t))$ satisfies $W(t)>0$ for all $t \in [0,T)$. The idea is to study the rate of change of the Lipschitz function $W(t)$ and show that it must decrease to zero in finite time. Notice that before the work of Perelman it was not even known whether an arbitrary  Ricci flow in a manifold like $S^3$ would develop a singularity in finite time. It is quite remarkable that this can be proved using minimal surface techniques.

Let $\Sigma$ be a branched minimal two-sphere in $(M^3,g(t_0))$. We denote by $N$ its unit normal vector field. The scalar curvature of a metric $g$ is denoted by $R_g$. We use $K_\Sigma$ to indicate the Gauss curvature of the induced metric on $\Sigma$, and $A$ to indicate its second fundamental form.
The rate of change of the area of a minimal surface under Ricci flow was first considered by Hamilton \cite{hamilton-survey}. We have
\begin{align*}
\frac{d}{dt}_{|t=t_0}{\rm area}_{g(t)}(\Sigma) = \frac12 \int_\Sigma {\rm tr}_\Sigma (-2Ric_{g(t_0)}) d\Sigma\\
=-\frac12 \int_\Sigma R_{g(t_0)}d\Sigma - \int_\Sigma K_\Sigma d\Sigma  -\frac12 \int_\Sigma |A|^2 d\Sigma\\
\leq-\frac12 \int_\Sigma R_{g(t_0)}d\Sigma - 4\pi - 2\pi \sum_i b_i,
\end{align*}
where we have used the Gauss equation and the Gauss-Bonnet Theorem with finitely many branch points $p_i$ of orders $b_i>0$.  Hence
$$
\frac{d}{dt}_{|t=t_0}{\rm area}_{g(t)}(\Sigma) \leq -4\pi -\frac12 \int_\Sigma R_{g(t_0)}d\Sigma.
$$

On the other hand, the scalar curvature of a Ricci flow $g(t)$ satisfies the evolution equation:
$$
\frac{\partial}{\partial t}R_{g(t)} = \Delta_{g(t)} R_{g(t)} + 2 |Ric_{g(t)}|^2.
$$
Therefore $\frac{\partial}{\partial t}R_{g(t)} \geq  \Delta_{g(t)} R_{g(t)} + \frac23 R_{g(t)}^2$, and the maximum principle for parabolic equations implies
$$
\min_M R_{g(t)} \geq \frac{1}{\frac{1}{\min_M R_{g(0)}}-\frac23 t}
$$
for all $t\in [0,T)$.

Putting things together, we get that 
\begin{equation*}
\frac{d}{dt}_{|t=t_0}{\rm area}_{g(t)}(\Sigma) \leq -4\pi +\frac{3}{4(t_0+C)} {\rm area}_{g(t_0)}(\Sigma)
\end{equation*}
for any branched minimal sphere in $(M,g(t_0))$, where the constant $C$ depends only on $g(0)$.

By comparison arguments and Theorem \ref{good.sweepouts}, Colding and Minicozzi proved:
\begin{theorem}[\cite{Colding-Minicozzi1}]\label{width.evolution.theorem}
Let $M^3$ be a closed orientable three-manifold that is prime and non-aspherical, and suppose that  $\varphi:S^2 \times [0,1]\rightarrow M$ induces a homotopically nontrivial map $\hat{\varphi}: S^3 \rightarrow M^3$. Then, for any Ricci flow $g(t)$ on $M$, we have
\begin{equation*}\label{width.derivative}
\frac{d}{dt}W(t)\leq -4\pi +\frac{3}{4(t+C)}W(t)
\end{equation*}
 in the sense of the limsup of forward difference quotients, where $W(t)=W(\varphi,g(t))$ and the constant $C>0$ depends only on $g(0)$.
\end{theorem}

This implies
$$
W(t)(t+C)^{-3/4} \leq -16\pi(t+C)^{1/4} + C'
$$
for some constant $C'$ independent of $t$. Since we always have $W(t)>0$, the flow cannot exist for all time.  This already implies that any Ricci flow in a prime non-aspherical orientable closed three-manifold must develop  a singularity in finite time. By arguing as in Perelman \cite{PERELMAN03B}, Theorem \ref{width.evolution.theorem} holds for Ricci flow with surgeries as well, thereby implying that any Ricci flow with surgery on one of these manifolds must become extinct in finite time.

\section{Conformally invariant variational problems}\label{willmore.mobius}

In this section we will describe some recent applications of the variational theory of minimal surfaces to global problems about the geometry of surfaces and links in three-space. 
We begin by describing  the min-max procedure used to produce the minimal hypersurface $\Sigma$ of Theorem \ref{pitts.theorem}. For simplicity, we restrict to the case $n=3$. 

A family of closed surfaces  $\{\Sigma(t)\}_{t \in [0,1]}$ of $M^3$  (where closed surface here means a two-dimensional integral cycle) is called a {\it sweepout} if we can write $\Sigma(t)=\partial \Omega(t)$ with the domains (integral 3-currents) $\Omega(t)$ varying continuously and satisfying $\Omega(0)=0$ and $\Omega(1)=M$. In particular we have $\Sigma(0)=\Sigma(1)=0\in \mathcal{Z}_2(M^3)$. A standard example of a sweepout is obtained by choosing a Morse function $f:M \rightarrow \mathbb{R}$, with  $f(M)=[0,1]$, and considering
\begin{align*}
 \Omega(t) = \{x \in M: f(x)<t\}   \\
 \Sigma(t)=\partial\Omega(t)=\{x \in M: f(x)=t\}.
 \end{align*}

We denote by  $\Pi_1$ the class of all sweepouts $\{\Sigma(t)\}_{t \in [0,1]}$ of $M^3$, and define the min-max
 invariant called the {\it width} of  $\Pi_1$:
$$
W(\Pi_1) = \inf_{\{\Sigma(t)\} \in \Pi_1} \sup_{t \in [0,1]} {\rm area}(\Sigma(t)).
$$
Since there will be some $t_0\in [0,1]$ such that ${\rm vol}(\Omega(t_0))={\rm vol}(M)/2$, the isoperimetric inequality implies that $W(\Pi_1)>0$. The minimal surface $\Sigma$ produced by Theorem \ref{pitts.theorem}, which could have several connected components  with integer multiplicities, is constructed  to satisfy 
$$
{\rm area}(\Sigma) = W(\Pi_1).
$$

If the ambient manifold is the three-dimensional unit sphere $S^3\subset \mathbb{R}^4$,  the minimal surface
produced by doing min-max over the class $\Pi_1$   is, modulo rotations, the equator or great sphere:
$\overline{\Sigma}=S^3 \cap \{x_4=0\}.$
In other words, $W(\Pi_1) = 4\pi$ and is achieved precisely by the great spheres.

 As a consequence, we have:
\begin{theorem}[$4\pi$ Theorem]
Let $\Phi:I \rightarrow \mathcal{Z}_2(S^3)$ be a  sweepout of $S^3$. Then there exists $y\in [0,1]$ such that
$${\rm area}(\Phi(y))\geq  4\pi.$$
\end{theorem}

The simplest minimal surface in $S^3$ after the equator is the Clifford torus $$\hat{\Sigma}=S^1(\frac{1}{\sqrt{2}}) \times S^1(\frac{1}{\sqrt{2}}),$$
with   area $2\pi^2$ and Morse   index 5. In  fact,  the Clifford torus can be characterized by its index:

\begin{theorem}[Urbano, \cite{urbano}]\label{urbano.theorem} Let $\Sigma \subset S^3$ be a smooth closed minimal surface of genus $g\geq 0$ and  ${\rm index}(\Sigma) \leq 5$. Then $\Sigma$ is either a great sphere (with index 1) or the Clifford torus (with index 5), up to ambient isometries.
\end{theorem}

The previous discussion implies that  the great sphere can appear as a min-max minimal surface, by considering one-parameter sweepouts. The question we posed ourselves, and that it turned out to be key to the solution of well-known global problems in conformal geometry such as the Willmore conjecture, was the following:

\medskip

{\bf Question:} {\em Is it possible to produce the Clifford torus  by min-max methods?}

\medskip

We have answered this question affirmatively by discovering a natural class of five-parameter families of surfaces in $S^3$ with interesting topological properties. The families we have discovered  are   parametrized by a map $\Phi$ defined on the 5-cube  $I^5$, and satisfy:
\begin{itemize}
\item[(A1)] $\Phi(x,0)=\Phi(x,1) = 0$ (trivial surface) for any $x\in I^4$,
\item[(A2)] $\{\Phi(x,t)\}_{t\in [0,1]}$ is the standard sweepout of $S^3$ by oriented round spheres centered at some $Q(x) \in S^3$, for any $x\in \partial I^4$,
\item[(A3)] $\Phi(x,1/2) = \partial B_{\pi/2}(Q(x))$, for any $x\in \partial I^4$,
\item[(A4)] ${\rm deg}(Q) \neq 0,$
\item[(A5)] there is no concentration of area:
$$\lim_{r\to 0}\sup \left\{{\rm area}(\Phi(x)\cap B_r(p)):\,p\in S^3, x \in I^5 \right\} =0.$$
\end{itemize}

The property (A4) above is saying that the restriction of the map $\Phi$ to  $\partial I^4 \times \{1/2\}$ is a homotopically nontrivial map
into  the space of oriented great spheres, which is homeomorphic to $S^3$.  This is the crucial topological condition that will rule out the possibility of
producing great spheres by min-max over families homotopic to $\Phi$.  The property (A5) is more of a technical nature. It is used in \cite{marques-neves} in the proofs of the interpolation statements.

The min-max theory developed jointly with Neves   in \cite{marques-neves} implies:

\begin{theorem}[$2\pi^2$ Theorem, \cite{marques-neves}] Let $\Phi:I^5 \rightarrow \mathcal{Z}_2(S^3)$ be a continuous map in the flat topology satisfying the properties (A1)-(A5) above. Then there must exist $y \in I^5$ with 
$${\rm area}(\Phi(y)) \geq 2\pi^2.$$
\end{theorem}

The surprising thing is that this is intimately related to the Willmore conjecture, a problem that we briefly discussed in Section \ref{recent.applications}. Let us describe this connection now.

First note that the torus $\Sigma_{\sqrt{2}}$ of Section \ref{recent.applications} is very special. In order to see this we need to recall that the stereographic projection  $\pi : S^3 \setminus \lbrace p \rbrace \rightarrow \mathbb{R}^3$, $p\in S^3$, is a conformal transformation. Since the Willmore energy is conformally invariant, we can  formulate the Willmore conjecture equivalently as a problem for surfaces in the three-sphere $S^3$.  If $\Sigma \subset S^3\setminus \lbrace p\rbrace$,  we can calculate the energy of the projection
 $\tilde{\Sigma} = \pi(\Sigma) \subset \mathbb{R}^3$:
  \begin{equation*}
\int_{\tilde{\Sigma}} \tilde{H}^2 d\tilde{\Sigma} = \int_{\Sigma} (1+H^2)\, d\Sigma,
\end{equation*}
where $H$ now denotes the mean curvature of $\Sigma$ with respect to the spherical geometry. 

Hence we define the {\it Willmore energy} of $\Sigma \subset S^3$ by
 $$\mathcal{W}(\Sigma) = \int_{\Sigma} (1+H^2)\, d\Sigma.$$

One advantage of considering the problem for surfaces in $S^3$ is that a relation with the area functional becomes apparent. It follows immediately from the above definition that for $\Sigma \subset S^3$ one always has $\mathcal{W}(\Sigma) \geq  {\rm area}(\Sigma)$, and $\mathcal{W}(\Sigma) = {\rm area}(\Sigma)$ if and only if $\Sigma$ is a minimal surface. If $\hat{\Sigma}$ denotes the Clifford torus  $S^1(1/\sqrt{2})\times S^1(1/\sqrt{2}) \subset S^3$, 
it is no coincidence that if we choose the right stereographic projection we will have $\pi(\hat{\Sigma})=\Sigma_{\sqrt{2}}$.

The Willmore conjecture has a long history of partial results. We refer the reader to our paper \cite{marques-neves} for an account.  We were motivated by the work of
Ros \cite{ros} on the antipodally symmetric case (see Topping \cite{topping} for a different proof), who used an area estimate we will mention in the sequel (Theorem \ref{area.estimate}).

A result of particular relevance to our approach is:
\begin{theorem}[Li and Yau, \cite{li-yau}] If $F:\Sigma \rightarrow S^3$ is an immersion and there exists $p\in S^3$ such that $\# F^{-1}(p)=k$, then $\mathcal{W}(\Sigma)\geq 4\pi k$.
In particular, if $\Sigma$ is not embedded then $\mathcal{W}(\Sigma)\geq 8\pi$.
\end{theorem}

Because of the result of Li and  Yau, we may assume the surface is embedded. Together with Neves we proved the following theorem that implies the Willmore conjecture:

\begin{theorem}[\cite{marques-neves}]\label{theorem.A}
Let $\Sigma\subset S^3$ be a closed, embedded smooth surface with  genus $g\geq 1$. Then $\mathcal{W}(\Sigma)\geq 2\pi^2$, and $\mathcal{W}(\Sigma)= 2\pi^2$ if and only if $\Sigma$ is a conformal image of the  Clifford torus.
\end{theorem}

\begin{remark}
The theorem above is completely equivalent to Theorem \ref{theorem.A.euclidean} mentioned in Section \ref{recent.applications}.
\end{remark}

The existence of a  torus that minimizes the Willmore energy was established by Simon \cite{simon93}. His work
was later extended to surfaces of higher genus by Bauer and Kuwert \cite{bauer-kuwert} (see also \cite{kusner96}). Very little is known about these surfaces and their energies. It is known that the minimum energy $\omega_g$ for orientable closed surfaces of genus $g$ in $\mathbb{R}^3$   is bigger than $2\pi^2$, less than $8\pi$, and converges to $8\pi$ as $g\rightarrow \infty$ \cite{kuwert-li-schatzle}.

It is also interesting to study general critical points of the Willmore energy, called Willmore surfaces. 
These are closed surfaces in $\mathbb{R}^3$ that satisfy the fourth-order Euler-Lagrange equation:
$$\Delta H + 2(H^2-K)H=0,$$
where $K$ denotes the Gauss curvature.  The simplest examples are stereographic projections of minimal surfaces of $S^3$, but there are many more.    Bryant \cite{bryant2} found and classified  all critical points of genus zero, and Pinkall \cite{pinkall} constructed infinitely many embedded Willmore tori  that are not the conformal image of  a minimal surface. The understanding of the analytical aspects of the Willmore equation has  been greatly improved in recent years thanks to the work of Kuwert-Sch\"atzle (e.g. \cite{kuwert-schatzle}) and Rivi\`ere (e.g. \cite{riviere}). The lecture notes of Rivi\`ere \cite{riviere-notes} provide a great introduction to the  analysis behind conformally invariant variational problems. 
We point also that the Willmore energy appears in many other interesting contexts, like in general relativity (as the main term in the definition of the Hawking mass) and in  relation with the renormalized area functional in the AdS/CFT correspondence, as recently studied by Alexakis and Mazzeo \cite{alexakis-mazzeo}.

In what follows we give an idea of the proofs of the $2\pi^2$ Theorem and Theorem \ref{theorem.A}, emphasizing their intimate relation.
We start  by constructing, for 
each embedded closed surface $\Sigma\subset S^3$, a {\it canonical family} of surfaces $\Sigma_{(v,t)}\subset S^3$, where $(v,t) \in B^4\times (-\pi,\pi)$, with the properties that:
\begin{itemize}
\item $\Sigma_{(0,0)}=\Sigma$,
\item  ${\rm area}(\Sigma_{(v,t)})\leq \mathcal{W}(\Sigma)$ for every $(v,t) \in B^4\times (-\pi,\pi)$.
\end{itemize}
Here $B^4\subset \mathbb{R}^4$ denotes the open unit ball. Each surface $\Sigma_{(v,t)}$ is an equidistant surface of some conformal image of $\Sigma$.

The fact that ${\rm area}(\Sigma_{(v,t)})\leq W(\Sigma)$ for every $(v,t)$ follows by combining the conformal invariance of the Willmore energy with the area estimate:
\begin{theorem}\label{area.estimate}
Let $\Sigma_t$, $t\in (-\pi,\pi)$, be an equidistant surface of an embedded closed surface $\Sigma \subset S^3$. Then
$$
{\rm area}(\Sigma_{t})\leq \mathcal{W}(\Sigma).
$$
Moreover, if $\Sigma$ is not a geodesic sphere and
\begin{equation*}
{\rm area }(\Sigma_{t})=\mathcal{W}(\Sigma),
\end{equation*}
then $t=0$ and $\Sigma$ is a minimal surface.
\end{theorem}

\begin{remark}
Theorem \ref{area.estimate} is a particular case of more general estimates proved by Heintze and Karcher 
\cite{heintze-karcher}. The statement above was used in connection with the Willmore problem  by  Ros \cite{ros}.
\end{remark}

The next step  is to understand the geometric and topological properties of the canonical family, especially the behavior as we let the parameter $(v,t)$ converge to the boundary of the parameter space $\partial \left(\overline{B}^4\times [-\pi,\pi]\right)$. We prove that:
\begin{itemize}
\item for any sequence $\{(v_i,t_i)\}$ that converges to  $\partial \left(\overline{B}^4\times [-\pi,\pi]\right)$, a subsequence of  $\{\Sigma_{(v_i,t_i)}\}$  converges in the flat topology to some round sphere (possibly trivial) in $S^3$.
\end{itemize}

The subtle case to consider is when $v$ converges to a point $p$ on the surface $\Sigma$, because then the limit is not unique and depends on the angle of convergence at which $v$ approaches $p$.  
We perform a blow-up procedure along the surface $\Sigma$ to solve this problem of  nonuniqueness of the limit. After reparametrizing the family, and observing that $\overline{B}^4 \times [-\pi,\pi]$ is homeomorphic to $I^5$, we get a family $\Phi:I^5 \rightarrow \mathcal{Z}_2(S^3)$ with properties (A1)-(A3)  and (A5) above, and an explicit center map $Q:\partial I^4 \rightarrow S^3$. 

The main topological ingredient is the discovery that the genus of the original surface $\Sigma$ can be read off the topological properties of the canonical family at the boundary:
\begin{theorem}[\cite{marques-neves}]
${\rm deg}(Q)= {\rm genus}(\Sigma)$.
\end{theorem}

Summarizing, we obtain:
\begin{theorem}[\cite{marques-neves}]\label{modified.family}
 Let $\Sigma \subset S^3$ be an embedded closed surface of genus $g\geq 1$. The map $\Phi:I^5 \rightarrow \mathcal{Z}_2(S^3)$ is continuous in the flat topology, satisfies the properties (A1)-(A5) and
$$
\sup\{{\rm area}(\Phi(x)):x\in I^5\}\leq \mathcal{W}(\Sigma).
$$
\end{theorem}

 Informally, the min-max family $\Phi$ can be thought of as an element of the relative homotopy group $\pi_5(\mathcal{S}, \mathcal{G})$, where
  $\mathcal{S}$ denotes the space of two-surfaces in $S^3$ and $\mathcal{G}$ denotes the space of round spheres.
  
  Given a smooth, embedded, closed surface $\Sigma \subset S^3$ with genus $g\geq 1$, it follows from the above properties that the
  map $\Phi:I^5 \rightarrow \mathcal{Z}_2(S^3)$ satisfies all the assumptions of the $2\pi^2$ Theorem. Therefore there exists $y \in I^5$ such that
  ${\rm area}(\Phi(y)) \geq 2\pi^2$. Since ${\rm area}(\Phi(x)) \leq \mathcal{W}(\Sigma)$ for every $x\in I^5$, we get $\mathcal{W}(\Sigma)\geq 2\pi^2$. In particular this proves the Willmore conjecture. The rigidity statement can be derived from the particular structure of the canonical family and the equality case in the area estimate of Ros.

We apply the Almgren-Pitts min-max theory for the area functional in order to prove the $2\pi^2$ Theorem. Given a family $\Phi:I^5 \rightarrow \mathcal{Z}_2(S^3)$ with properties (A1)-(A5), we consider $\Pi$ the homotopy class of $\Phi$ relative to $\partial I^5$. 
Of course we have
$$
\sup\lbrace {\rm area}(\Phi(x)): x\in \partial I^5\rbrace = 4\pi.
$$
We first rule out great spheres as possible min-max surfaces for $\Pi$, by proving:

\begin{theorem}[\cite{marques-neves}]\label{no.great.spheres}
$L(\Pi)>4\pi$.
\end{theorem}

The proof is topological and goes by contradiction, by assuming $L(\Pi)=4\pi$. In order to illustrate the idea,  we suppose there is a map $\tilde{\Phi}\in \Pi$ with 
\begin{equation*}
\sup_{x \in I^5} {\rm area}(\tilde{\Phi}(x)) = 4\pi.
\end{equation*}
In particular, for any given continuous path $\gamma:[0,1]\rightarrow I^5$ with $\gamma(0)\in I^4\times \{0\}$ and $\gamma(1)\in I^4 \times \{1\}$,  $\{\tilde{\Phi} \circ \gamma\}$ is optimal as a one-parameter sweepout of $S^3$. Therefore it must contain a great sphere. One could argue then that there must exist a 4-dimensional submanifold $R \subset I^5$, separating the top from the bottom of $I^5$, such that
\begin{itemize}
\item $\tilde{\Phi}(y)$ is a great sphere for any $y\in R$, 
\item $\partial R= \partial I^4 \times \{1/2\}$.
\end{itemize}
If $\tilde{Q}(y)$ denotes the center of the great sphere $\tilde{\Phi}(y)$, for $y\in R$, and since $\tilde{\Phi}=\Phi$ on $\partial I^5$, we get
$$
[Q_\#(\partial I^4)]=[\tilde{Q}_\#(\partial R)]=[\partial \tilde{Q}_\#(R)]=0\in H_3(S^3,\mathbb{Z}).
$$
But $Q_\#(\partial I^4)={\rm deg}(Q) \cdot S^3$, and we reach a contradiction since ${\rm deg}(Q)\neq 0$. The complete proof is more involved and can be found in \cite{marques-neves}.

Once we know that $L(\Pi)>4\pi$, here is how we prove  that in reality $L(\Pi) \geq 2\pi^2$. By applying min-max theory to $\Pi$, we get the existence of a smooth embedded minimal surface $\Sigma^\prime\subset S^3$ such that
\begin{equation*}
L(\Pi)={\rm area}(\Sigma^{\prime}) >4\pi.
\end{equation*}
If $\Sigma^\prime$ has multiplicity bigger than one, then $L(\Pi)\geq 8\pi$. If the multiplicity is one, and since by Almgren \cite{almgren66} the only minimal spheres in $S^3$ are the great spheres, we get ${\rm genus}(\Sigma^\prime)\geq 1$.

The result follows once we show that the nonspherical minimal surface of lowest area $\hat{\Sigma}$ in $S^3$ is the Clifford torus. The idea is to use Urbano's theorem. If we had ${\rm index}(\hat{\Sigma})\geq 6$, we would be able to slightly  perturb its canonical family and produce, by min-max theory and Theorem \ref{no.great.spheres}, a minimal surface whose area is strictly between $4\pi$ and ${\rm area}(\hat{\Sigma})$. Contradiction, hence ${\rm index}(\hat{\Sigma})\leq 5$ and $\hat{\Sigma}$ must be the Clifford torus by Urbano's theorem. Of course, ${\rm area}(\Sigma') \geq {\rm area}(\hat{\Sigma})$. This finishes a sketch of the proof that $L(\Pi) \geq 2\pi^2$, and this implies the $2\pi^2$ Theorem.

\begin{remark}
Note that in the process of proving the $2\pi^2$ Theorem we have showed Theorem \ref{theorem.B}. 
\end{remark}

\subsection{Links}
In order to prove Theorem \ref{fhw.conjecture.theorem}, we will again use the $2\pi^2$ Theorem. The basic observation is that if $g$ denotes the Gauss map of a link $(\gamma_1,\gamma_2)$ contained in $S^3$, i.e., the map  $g:S^1 \times S^1 \rightarrow S^3$ defined by
$$
g(s,t) = \frac{\gamma_1(s)-\gamma_2(t)}{|\gamma_1(s)-\gamma_2(t)|},
$$
then 
$
|{\rm Jac} \, g|(s,t) \leq \frac{|\gamma_1'(s)||\gamma_2'(t)|}{|\gamma_1(s)-\gamma_2(t)|^2}.
$
Hence
$
{\rm area}(g(S^1\times S^1)) \leq E(\gamma_1,\gamma_2).
$

By applying conformal transformations to the curves $\gamma_1$, $\gamma_2$ of a link $(\gamma_1,\gamma_2)$ in $S^3 \subset \mathbb{R}^4$, and considering the associated Gauss maps, we get a 5-parameter family of surfaces (parametrized tori) in $S^3$. The key again is to analyze the boundary behavior. After an extension, we get a family with the same basic properties of the canonical family for the Willmore problem, and such that the area  of any surface in the family is bounded above by the M\"{o}bius energy of the link.  We prove that if $|{\rm lk}(\gamma_1,\gamma_2)| = 1$ then  the center map $Q:\partial I^4 \rightarrow S^3$ associated with the family satisfies $|{\rm deg}(Q)|=1$. Therefore the $2\pi^2$ Theorem applies and we conclude the existence of at least one surface in the family with area greater than or equal to $2\pi^2$. This establishes the inequality $E(\gamma_1,\gamma_2) \geq 2\pi^2$, and after some extra work one can also prove the rigidity part. 

\section{Further directions}\label{further}

In this section we describe some recent advances in min-max theory and discuss some future directions. We start with the problem of counting minimal hypersurfaces 
in Riemannian manifolds. 
We will give an informal overview of the proof of Theorem \ref{infinitely.theorem}. 

Let $M^{n+1}$ be a Riemannian manifold as in Theorem \ref{infinitely.theorem}. The homotopy groups of the space of modulo $2$  $n$-cycles in $M$, $\mathcal{Z}_{n}(M,\mathbb{Z}_2)$, can be computed through the work of Almgren \cite{almgren}. All homotopy groups vanish but the first one: $\pi_1(\mathcal Z_{n}(M,\mathbb{Z}_2)) = \mathbb{Z}_2$, just like for the topological space  $\mathbb{RP}^{\infty}$. Let $\bar \lambda \in H^1(\mathcal Z_{n}(M,\mathbb{Z}_2), \mathbb{Z}_2)$ be the generator.  


Gromov \cite{gromov0, gromov, gromov2} and Guth \cite{guth} have studied  continuous maps $\Phi$ from a simplicial complex $X$ into $\mathcal Z_{n}(M,\mathbb{Z}_2)$ that detect $\bar \lambda^p$, in the sense that $\Phi^{*}(\bar \lambda^p)\neq 0$.  Here $\bar \lambda^p$ denotes the $p$-th cup power of $\bar \lambda$. We call these maps {\it $p$-sweepouts}. An example can be given by starting with a Morse function $f: M \rightarrow \mathbb{R}$. The open set  $\{x\in M: f(x)<t\}$ has finite perimeter for all $t$, hence we have a well-defined element $$f^{-1}(t)=\partial \{x\in M: f(x)<t\}\in \mathcal Z_n(M;\mathbb{Z}_2).$$ 
For each  $a=(a_0,\ldots,a_p) \in \mathbb{R}^{p+1}$, $|a|=1$, we consider the polynomial $P_a(t)=\sum_{i=0}^p a_i t^i$ and  define the map  
$$\Psi:\{a\in \mathbb{R}^{p+1}:|a|=1\}\rightarrow \mathcal{Z}_n(M;\mathbb{Z}_2)$$
by
 $$
\Psi(a_0,\ldots,a_p)=\partial\left \{x\in M:P_a(f(x))<0\right\}.$$
The fact that we are using  $\mathbb{Z}_2$ coefficients implies that  $\Psi(a)=\Psi(-a)$, and therefore $\Psi$ induces a  map $\Phi:\mathbb{RP}^p\rightarrow \mathcal{Z}_n(M;\mathbb{Z}_2).$ It satisfies  $\Phi^{*}(\bar \lambda^p)\neq 0$.

We denote by  $\mathcal P_p$ the space of all maps  that detect $\bar\lambda^p$, and define the min-max invariant:
 \begin{equation*}\label{sublinear}
 \omega_p(M):=\inf_{\Phi\in \mathcal P_p}\sup_{x\in {\rm dmn}(\Phi)} {\rm area}(\Phi(x)),
 \end{equation*}
 where ${\rm dmn}(\Phi)$ stands for the domain of $\Phi$.
 The asymptotic behavior of the min-max $n$-volumes $\omega_p(M)$ as $p \rightarrow \infty$ has been studied previously by Gromov and Guth. The following result was proven by Gromov in \cite[Section 4.2.B]{gromov0}, and by Guth in \cite{guth} via an elegant bend--and--cancel argument.
 
\begin{theorem}\label{upper.bound}   There exists a constant $C=C(M)>0$ so that $$\omega_p(M)\leq  C p^{\frac{1}{n+1}}$$
for every $p\in \mathbb{N}$.
\end{theorem}

{\it Remark:} The lower bound $  \omega_p(M)\geq  C' p^{\frac{1}{n+1}}$ also holds for some constant $C'=C'(M)>0$ (Gromov \cite{gromov}, Guth \cite{guth}).
 
 \medskip
 
 We use Lusternik-Schnirelmann theory to show that if $\omega_p=\omega_{p+1}$ then there are infinitely many embedded minimal hypersurfaces. This is inspired by the particular structure of the cohomology ring of a projective space.  Details can be found in our paper \cite{marques-neves-infinitely}.
 
 We then prove Theorem \ref{infinitely.theorem} by contradiction,  assuming that the set $\mathcal{L}$ of all  smooth, closed, embedded minimal hypersurfaces in $M$ is finite. This implies that the sequence  $\{\omega_p(M)\}_{p\in \mathbb{N}}$ is strictly increasing. Since the support of the Almgren-Pitts min-max minimal surface is always embedded, the Frankel property implies that it must have the form $k \cdot \Sigma$ for some $\Sigma \in \mathcal{L}$. We conclude, by applying min-max theory to the classes of $p$-sweepouts and Theorem \ref{upper.bound}, that 
 \begin{align*}
\#\{a=k |\Sigma|: k\in \mathbb{N},\, \Sigma \in \mathcal{L},\, k |\Sigma| \leq Cp^{\frac{1}{n+1}}\}\\
\geq \#\{\omega_k(M):k=1,\ldots, p\}= p.
\end{align*}
 Since $\mathcal{L}$ is assumed to be finite, we get a contradiction for sufficiently large $p$.

\subsection{Open problems}We describe some open problems related to the variational theory of minimal surfaces.
Some of these questions are well-known and others arose from extensive discussions with Andr\'{e} Neves.

The precise relation between  the Morse index of the min-max minimal surface and the number of parameters is not known  in the Almgren-Pitts min-max theory.  In general, one should expect that ${\rm index}(\Sigma) \leq k$, where $k$ is the number of parameters. This is a subtle question, specially because  of the phenomenon of multiplicity. In \cite{zhou}, Zhou proves this for $k=1$ in the case of compact manifolds $M^n$ of positive Ricci curvature, with $3\leq n\leq 7$. This extends the results for $n=3$ of the author and Neves \cite{marques-neves-duke}.

 It is also natural to ask whether one can control the topology of $\Sigma$ produced by Theorem \ref{pitts.theorem}. This problem has been studied in the three-dimensional case through  variants of the Almgren-Pitts theory due to Simon-Smith \cite{smyth} and Colding and De Lellis \cite{colding-delellis}, in which the one-parameter sweepouts are made of actual smooth closed surfaces (perhaps up to finitely many singularities) with bounded genus. Simon and Smith proved in \cite{smyth}, by considering sweepouts of the form $\Sigma(t)=\psi(t,\overline{\Sigma}(t))$, where $\overline{\Sigma}(t)=\{x\in S^3: x_4=2t-1\}$ is the standard foliation by round spheres and $\psi: [0,1] \times S^3 \rightarrow S^3$ is an ambient isotopy, that every Riemannian metric on $S^3$ admits a minimal embedded two-sphere. 

These theories can also deal with the case of sweepouts induced by Heegaard splittings.  In the 1980s Pitts and Rubinstein  made a number of conjectures about the index and genus of the min-max minimal surface obtained this way. They have conjectured, for instance, that after performing finitely many surgeries of controlled type to the Heegaard surfaces, each remaing component should be isotopic to a component of the min-max minimal surface or to a double cover of a component. A proof of this last conjecture, which in turn implies basically optimal genus bounds, has been recently presented by Ketover \cite{ketover} (a previous result on genus bounds can be found in De Lellis and Pellandini \cite{delellis-genus}). It should be interesting  to bound the Betti numbers of a min-max minimal hypersurface in higher dimensions. A shorter proof of Theorem \ref{pitts.theorem}, still in the setting of one-parameter sweepouts, was given by De Lellis and Tasnady in \cite{delellis-tasnady}. 

Are the Clifford hypersurfaces 
$$\Sigma_{k+1,l+1}=S^{k}(\sqrt{\frac{k}{k+l}})\times S^{l}(\sqrt{\frac{l}{k+l}})\subset S^{n}(1),$$
$k+l=n-1$, the only non-equatorial minimal hypersurfaces of $S^{n}$ with Morse index less than or equal to $(n+2)$?
If $n=3$ this follows from Urbano \cite{urbano} (Theorem \ref{urbano.theorem}). See Perdomo \cite{perdomo} for the antipodally symmetric case.

The minimal hypersurfaces $\Sigma_{m,m}$ and $\Sigma_{m,m+1}$  are conjectured (Solomon) to have least area among all non-equatorial minimal hypersurfaces of $S^{2m-1}$ and $S^{2m}$, respectively. Similarly, the respective cones $C_{m,m}$ and $C_{m,m+1}$ over $\Sigma_{m,m}$ and $\Sigma_{m+1,m}$ should be the nontrivial minimal hypercones of least possible density. For $m=1$ these conjectures follow from the results of \cite{marques-neves}. See also Ilmanen and White \cite{ilmanen-white} for some partial results. 
 
The general behavior of singularities of minimal submanifolds is very little understood.  For instance, we do not know whether  there exists a properly embedded minimal surface in $B^3_1(0)\setminus \{0\}$ that does not extend smoothly to the origin.  Another question (Schoen) is whether singularities of area-minimizing hypersurfaces are stable under small perturbations of the data. See N. Smale \cite{smale} for the case $n=8$.  
 
 Gromov \cite{gromov0} has proposed to consider the sequence $\{\omega_p(M)\}_{p\in\mathbb{N}}$ as a nonlinear analogue of the Laplace spectrum of $M$. In particular one can ask (\cite[Section 8]{gromov}, \cite[Section 5.2]{gromov2}) whether a Weyl Law holds for the sequence of numbers 
 $\{\omega_p(M)\}_{p\in\mathbb{N}}.$ More precisely, if
\begin{equation*}\label{weyl}
\lim_{p\to\infty}\omega_p(M)p^{-\frac {1}{n+1}}=a(n)({\rm vol}(M,g))^{\frac{n}{n+1}},
\end{equation*}
where  $a(n)$ is a constant that depends only on $n$.

In \cite{marques-neves-infinitely}, we put this analogy forward by considering sweepouts whose surfaces are zero sets of linear combinations  of eigenfunctions. Note that a conjecture of Yau \cite{yau1} states that $$c^{-1}\sqrt \lambda_p\leq  \mathcal H^n(\{\phi_p=0\})\leq c\sqrt \lambda_p,$$ 
where $c=c(M,g)>0$, $\phi_p$ is the $p$-th eigenfunction and $\lambda_p$ is the $p$-th eigenvalue of the Laplacian.
We conjecture  that generically the minimal surfaces given by Theorem \ref{infinitely.theorem} should form a sequence $\{\Sigma_p\}_{p\geq 1}$, where $\Sigma_p$ has index $p$, multiplicity one and area going to infinity. The first Betti number and the index should be linearly related, and by analogy with nodal sets of eigenfunctions, we conjecture that these minimal surfaces  should become equidistributed in space.

Finally, we mention that there are many interesting questions in the higher codimension case, where the second variation formula becomes a complicated object.  It would be interesting  to construct calibrated submanifolds by variational methods (see Schoen and Wolfson \cite{schoen-wolfson}), and to understand better the influence of the ambient curvature.
The only stable integral cycles of complex projective space are the formal sums of algebraic varieties (Lawson and Simons \cite{lawson-simons}).
   



\end{document}